\documentclass[a4paper,12pt]{amsart}
\usepackage{amsmath}
\usepackage{amsfonts}
\usepackage{amssymb}
\usepackage{amsthm}
\usepackage{multicol}
\newtheorem{thm}{Theorem}[section]
\newtheorem{cor}[thm]{Corollary}
\newtheorem{lemma}[thm]{Lemma}
\newtheorem{prop}[thm]{Proposition}
\theoremstyle{definition}
\newtheorem{remark}[thm]{Remark}
\numberwithin{equation}{section}

\def\al{\alpha}
\def\be{\beta}
\def\de{\delta}
\def\ep{\varepsilon}
\def\la{\lambda}
\def\si{\sigma}
\def\vp{\varphi}
\def\Z{\mathbb{Z}}

\def\R{\mathbb{R}}
\def\C{\mathbb{C}}
\def\N{\mathbb{N}}
\def\F{\mathcal{F}}

\newcommand{\rphis}[5]{\,_{#1}\vp_{#2} \left( \genfrac{.}{.}{0pt}{}{#3}{#4}
\ ;q, #5 \right)}
\newcommand{\tfo}[3]{\,_2\varphi_1 \left( \genfrac{.}{.}{0pt}{}{#1}{#2}
\ ;q, #3 \right)}
\newcommand{\ra}{\rightarrow}
\newcommand{\da}{\downarrow}
\newcommand{\usum}{\sum_{n=0}^\infty}

\begin{document}
\title[Al-Salam--Chihara polynomials]{
Self-adjoint difference operators and \\ symmetric Al-Salam--Chihara
polynomials}

\author{Jacob S.~Christiansen and Erik Koelink}
\address{California Institute of Technology, Mathematics 253-37, Pasadena, CA 91125, U.S.A.}
\email{stordal@caltech.edu}
\address{Technische Universiteit Delft, DIAM, PO Box 5031,
2600 GA Delft, the Netherlands}
\email{h.t.koelink@tudelft.nl}

\begin{abstract} The symmetric Al-Salam--Chihara polynomials
for $q>1$ are associated with an indeterminate moment problem. There
is a self-adjoint second order difference operator on $\ell^2(\Z)$
to which these polynomials are eigenfunctions. We determine the
spectral decomposition of this self-adjoint operator. This leads to
a class of discrete orthogonality measures, which have been obtained
previously by Christiansen and Ismail using a different method, and
we give an explicit orthogonal basis for the corresponding weighted
$\ell^2$-space. In particular, the orthocomplement of the
polynomials is described explicitly. Taking a limit we obtain all
the $N$-extremal solutions to the $q^{-1}$-Hermite moment problem, a
result originally obtained by Ismail and Masson in a different way.
Some applications of the results are discussed.
\end{abstract}

\maketitle

\begin{small}
\end{small}

\smallskip

\begin{small}
\end{small}


\section{Introduction}\label{sec:introduction}

The Al-Salam--Chihara polynomials were introduced originally by
Al-Salam and Chihara \cite{AlSaC} in their search for polynomials
satisfying a convolution identity, see \cite[\S 8]{AlSa} and also
\cite[Thm.~4.5]{KoelVdJ} for a generalization of the convolution
property of the Al-Salam--Chihara polynomials. They have been
studied by Askey and Ismail \cite[\S 3]{AskeI}, and it turns out
that they fit in the Askey-scheme of basic hypergeometric orthogonal
polynomials, see \cite{KoekS}. In case the base $q$ is $>1$ the
Al-Salam--Chihara polynomials may be related to an indeterminate
moment problem. The precise conditions on the parameters for this to
happen is given in \cite[Thm. 3.2]{AskeI}. When $q>1$ and the moment
problem is determinate, the orthogonality measure is purely discrete
\cite[(3.82)]{AskeI} and the dual polynomials are little $q$-Jacobi
polynomials. This observation can be found at various places in the
literature, see Groenevelt \cite[Rmk. 3.1]{Groe}, Atakishiyev and
Klimyk \cite[\S 3]{AtakK} and Rosengren \cite{Rose}.

The indeterminate moment problems arising in the Askey-scheme have
been classified and studied in \cite{ChriPhD}. In this paper we
study the symmetric Al-Salam--Chihara polynomials for $q>1$. They
are associated with an indeterminate moment problem, see
\cite{ChriI} for more information on various solutions of the moment
problem in question. In fact, they form a one-parameter extension of
the continuous $q^{-1}$-Hermite polynomials, which correspond to one
of the very few examples, if not the only example, of an
indeterminate moment problem for which all $N$-extremal measures are
known explicitly, see Ismail and Masson \cite{IsmaM} for details.
The main idea of this paper is to exploit the fact that the
Al-Salam--Chihara polynomials occur in the Askey-scheme. This
implies that they are eigenfunctions to a specific second order
difference operator $L$. We look for a suitable Hilbert space
${\mathcal H}$ of functions such that $L$ is a self-adjoint operator
on ${\mathcal H}$. The spectral decomposition of $L$ then gives
information on denseness of polynomials in ${\mathcal H}$, or on the
orthocomplement of the space of polynomials in ${\mathcal H}$. This
has been exploited in our previous paper \cite{ChriK} for the
Stieltjes--Wigert polynomials, but the case here is much simpler.
See also \cite{CiccKK} for the $q$-Laguerre polynomials and
\cite{KoelS} for the continuous dual $q^{-1}$-Hahn polynomials,
where in both cases the result occurred as a by-product of the study
of a certain self-adjoint operator.

In Section \ref{sec:ASCpolynomials} we introduce the symmetric
Al-Salam--Chihara polynomials, and we recall some of their
properties. Section \ref{sec:selfadjointdiffop} is the main part of
the paper. We first rewrite the second order difference operator $L$
as a self-adjoint operator acting on $\ell^2(\Z)$. We show that this
operator is trace-class, and even in every Schatten class, by
estimating its singular values. Next we list several solutions to
the corresponding difference equation using known contiguous
relations for basic hypergeometric series, see Lemma
\ref{lem:manysolutions}. With this information at hand we can give
the spectral decomposition for $L$ explicitly, using standard
techniques as described e.g. in \cite{Koel}. It turns out that the
spectrum has a positive discrete part corresponding to the symmetric
Al-Salam--Chihara polynomials with parameter $\be$, and a negative
discrete part corresponding to the symmetric Al-Salam--Chihara
polynomials with parameter $1/q^2\be$. The special cases $\be\da 0$
and $\be=1/q$ are related to the continuous $q^{-1}$-Hermite
polynomials. As explained in Section 4, taking the limit $\be\da 0$
we obtain the $N$-extremal solutions to the $q^{-1}$-Hermite moment
problem, a result originally due to Ismail and Masson \cite{IsmaM}.
The proof of Ismail and Masson is based on explicit descriptions of
the entire functions in the Nevanlinna parametrization and heavily
relies on theta function identities. We emphasize that our proof is
different, but that the outcome of the $N$-extremal measures is a
lucky coincidence. The case $\be=1/q$ is related to orthogonality
measures studied by Christiansen and Ismail \cite{ChriI}. We prove
some of their results in a different way, and we extend some of
their results as well. In particular, in \cite{ChriI} the derivation
of the measures $\la_\al^{(\be)}$ is based on the use of Bailey's
${}_6\psi_6$-summation. In the present setting, we can reverse the
argument to obtain special cases of Bailey's ${}_6\psi_6$-summation.


\section{Al-Salam--Chihara polynomials}\label{sec:ASCpolynomials}

In this section we fix the notation for the Al-Salam--Chihara
polynomials, and we recall some results for these polynomials. We
use the standard notation for basic hypergeometric series as in
Gasper and Rahman \cite{GaspR}. By switching to base $q^{-1}$, we
now assume that $0<q<1$. The Al-Salam--Chihara polynomials in base
$q^{-1}$ are, see e.g. \cite{AskeI}, \cite{KoekS},
\begin{align}
\label{eq:defASCpols} P_n(u;a,b \,|\, q^{-1}) &=
\frac{(ab;q^{-1})_n}{a^n} \,_{3}\vp_{2} \left(
\genfrac{.}{.}{0pt}{}{q^n, au, a/u}{ab, 0} \ ;q^{-1}, q^{-1}
\right) \\
\notag &= \frac{(-b)^n(1/ab;q)_n}{q^{\binom{n}{2}}}
\rphis{3}{1}{q^{-n},\, 1/au,\, u/a}{1/ab}
{\frac{aq^n}{b}} \\
\notag &= \frac{(-abu)^n(1/ab;q)_n}{q^{\binom{n}{2}}}
\rphis{3}{2}{q^{-n},\, 1/au,\, 1/bu} {1/ab,\, 0}{q} \\
\notag &= \frac{(-b)^n(u/b;q)_n}{q^{\binom{n}{2}}}
\rphis{2}{1}{q^{-n},\, 1/au}{bq^{1-n}/u}{\frac{aq}{u}},
\end{align}
which is a polynomial of degree $n$ in $\frac12 (u+u^{-1})$. In
\eqref{eq:defASCpols} we use \cite[Exercise 1.4(i)]{GaspR} to switch
from base $q^{-1}$ to base $q$ and \cite[(III.7), (III.8)]{GaspR} to
rewrite the result as a ${}_3\vp_2$- and ${}_2\vp_1$-series. Since
the Al-Salam--Chihara polynomials are symmetric in $a\leftrightarrow
b$ and $u \leftrightarrow u^{-1}$, we get several more expressions
for them.

The Al-Salam--Chihara polynomials are eigenfunctions to a second
order difference operator, which is a special case of the difference
operator for the Askey-Wilson polynomials by setting two parameters
equal to zero, see \cite{GaspR}, \cite[\S 3.8]{KoekS}. Explicitly,
$P_n(u)$ as in \eqref{eq:defASCpols} satisfies
\begin{equation}
\label{eq:diffeqASCpols} (q^n-1)\, P_n(u) =
A(u)\bigl[P_n(u/q)-P_n(u) \bigr] + A(u^{-1})\bigl[ P_n(uq)-P_n(u)
\bigr]
\end{equation}
with
\[
A(u) = \frac{(1-au)(1-bu)}{(1-u^2)(1-u^2/q)}.
\]
The three-term recursion generating $P_n(u)$ is, see \cite{AlSaC},
\cite[\S 3.8]{AskeI}, \cite[\S 3.8]{KoekS},
\begin{equation}
\label{eq:ttrASCpols}
(u+u^{-1})\, P_n(u) = P_{n+1}(u) + (a+b)q^{-n}\, P_n(u) +
(1- q^{-n})(1-abq^{1-n})\, P_{n-1}(u),
\end{equation}
with initial values $P_{-1}(u)=0$ and $P_0(u)=1$ .

In this paper we only deal with the symmetric Al-Salam--Chihara
polynomials, i.e. the case $a=-b$, although we comment on the
general case in Remark \ref{rmk:othercases}. By putting $u=i e^{-y}$
so that $\frac12 (u+u^{-1})= -i\sinh y$, we define
\begin{equation}
\label{eq:defsymASCpols}
Q_n(y):=Q_n(y;\be \,|\,q)= i^n P_n(ie^{-y};
\sqrt{\be},-\sqrt{\be}\,|\,q^{-1})
\end{equation}
for $\be>0$. The limiting case $\be\to 0$ can be obtained from the
explicit expression
\[
Q_n(y)=i^n\be^{n/2}\frac{(-ie^{-y}/\sqrt{\be};q)_n}{q^{\binom{n}{2}}}
\rphis{2}{1}{q^{-n},\,
-ie^y/\sqrt{\be}}{ie^yq^{1-n}\sqrt{\be}}{-iqe^y\sqrt{\be}}
\]
or by considering the limit transition $\be\to 0$ in the three-term
recurrence relation
\begin{equation*}
\label{eq:ttrsymmASCpols} (e^y-e^{-y})\,Q_n(y)=
Q_{n+1}(y)+q^{-n}(1-q^n)(1+\be q^{1-n})\,Q_{n-1}(y).
\end{equation*}
We introduce the symmetric Al-Salam--Chihara polynomials as
\begin{equation}\label{eq:defsymmhASCpols}
h_n^{(\be)}(x|q)=Q_n(y;\be\,|\,q), \quad x=\sinh y.
\end{equation}
They satisfy the recurrence relation
\begin{equation}\label{eq:ttrsymmhASCpols}
2x \, h_n^{(\be)}(x|q) = \, h_{n+1}^{(\be)}(x|q) + q^{-n}
(1-q^n)(1+\be q^{1-n})\, h_{n-1}^{(\be)}(x|q),
\end{equation}
and comparing this with the three-term recursion for the continuous
$q^{-1}$-Hermite polynomials, see \cite[(1.16)]{IsmaM},
\cite[(1.1)]{ChriI},
\begin{equation}\label{eq:ttrcontqinvHermitepols}
2x\, h_n(x|q) = h_{n+1}(x|q)\,  + \, q^{-n}(1-q^n)\, h_{n-1}(x|q),
\end{equation}
we see that $\lim_{\be\downarrow 0} h_n^{(\be)}(x|q) = h_n(x|q)$.
Moreover, we obtain $h_n^{(1/q)}(x|q) = h_n(x|q^2)$ by comparing
\eqref{eq:ttrsymmhASCpols} and \eqref{eq:ttrcontqinvHermitepols},
cf. \cite[\S 2]{ChriI}. Note that \eqref{eq:ttrsymmhASCpols} and
Favard's theorem, see e.g. \cite[\S 2]{AskeI}, \cite[\S 7.1]{GaspR},
\cite{Koel}, implies that the symmetric Al-Salam--Chihara
polynomials are orthogonal on $\R$. Moreover, by \cite[Thm.
3.2]{AskeI} it follows that the corresponding moment problem is
indeterminate.

Combining \eqref{eq:defsymASCpols} and \eqref{eq:diffeqASCpols} we
see that
\begin{equation}\label{eq:diffeqsymASCpols}
\begin{split}
q^n Q_n(y)=&\frac{1+\be e^{-2y}}{(1+e^{-2y})(1+e^{-2y}/q)}
\bigl(Q_n(y+\ln q)-Q_n(y)\bigr) \\ &\quad +Q_n(y)+ \frac{1+\be
e^{2y}}{(1+e^{2y})(1+e^{2y}/q)} \bigl(Q_n(y-\ln q)-Q_n(y)\bigr),
\end{split}
\end{equation}
and the corresponding difference operator will be the topic of
Section \ref{sec:selfadjointdiffop}.


\section{Spectral decomposition of a self-adjoint difference operator}
\label{sec:selfadjointdiffop}

In this section we introduce the self-adjoint operator $L$ that can
be obtained from \eqref{eq:diffeqsymASCpols} and start out studying
its general properties. Then we look for sufficiently many
eigenfunctions in order to determine the spectral decomposition
explicitly. As a by-product we obtain orthogonality relations
involving the symmetric Al-Salam--Chihara polynomials.

\subsection{Self-adjoint difference operator}
\label{ssec:selfadjointdiffop}

Fix $\al\in \R\backslash\{0\}$ and define
\begin{equation}\label{eq:defgridpoints}
x_k(\al) = \frac12 \Bigl(  \frac{1}{\al q^k} - \al q^k\Bigr), \qquad
k\in \Z.
\end{equation}
Note that $x_k(-\al)= x_{-k}(1/\al)$ and that $x_k(\al)\to\pm\infty$
as $k\to\pm\infty$. We denote by $Z(\al)$ the grid $\{x_k(\al)\mid
k\in\Z\}$ and have $Z(\al)=Z(\gamma)$ if and only if ${\al}/{\gamma}
\in q^\Z$ or $-\alpha\gamma\in q^\Z$. So it suffices to consider
$\al\in(q,1]$, and $\R=\cup_{\al\in (q,1]} Z(\al)$ as a union of
disjoint sets. Restricting to $\al\in (q,1]$ we see that
$Z(\al)=-Z(\al)$ only for $\al=\sqrt{q}$ and $\al=1$.

With the above notation the eigenvalue difference equation
\eqref{eq:diffeqsymASCpols} for $Q_n$ can be regarded as an operator
acting on $\F(\al):=\{ f\colon Z(\al)\to \C\}$. Defining $\de_l\in
\F(\al)$ by $\de_l\bigl(x_k(\al)\bigr)=\de_{k,l}$, we get a basis
for the finitely supported functions in $\F(\al)$ and the operator
$L=L(\al,\be)$ is given by
\begin{equation}\label{eq:defLalbeonFal}
\begin{split}
L\, \de_l = &\, \frac{1+\be\al^2 q^{2l+2}}{(1+\al^2
q^{2l+2})(1+\al^2 q^{2l+1})}\, \de_{l+1} \\ &\, + \left(1 -
\frac{1+\be\al^2 q^{2l}}{(1+\al^2 q^{2l})(1+\al^2 q^{2l-1})} -
\frac{1+\be \al^{-2} q^{-2l}}{(1+
\al^{-2}q^{-2l})(1+\al^{-2}q^{-2l-1})} \right) \de_l \\ &\, +
\frac{1+\be \al^{-2} q^{-2l+2}}{(1+
\al^{-2}q^{-2l+2})(1+\al^{-2}q^{-2l+1})}\, \de_{l-1}.
\end{split}
\end{equation}
To see how \eqref{eq:defLalbeonFal} is motivated from
\eqref{eq:diffeqsymASCpols}, restrict \eqref{eq:diffeqsymASCpols} to
the grid $Z(\al)$, i.e. replace $e^{-y}$ by $\al q^l$. Then identify
the polynomial $Q_n(y)$ with the element $\sum_l
Q_n\bigl(x_l(\al)\bigr)\, \de_l \in {\mathcal F}(\al)$, so that
$Q_n(y\pm \ln q)$ is identified with $\sum_l Q_n\bigl(x_{l\mp
1}(\al)\bigr)\, \de_l \in {\mathcal F}(\al)$. Inserting this in
\eqref{eq:diffeqsymASCpols} and shifting $l$ such that all terms
involve $Q_n\bigl(x_l(\al)\bigr)$, we see that the right-hand side
of \eqref{eq:diffeqsymASCpols} leads to the definition of $L$ in
\eqref{eq:defLalbeonFal}.

Note that \eqref{eq:diffeqsymASCpols} shows that the
Al-Salam--Chihara polynomials are eigenfunctions to $L$
corresponding to the eigenvalues $q^n$, $n\geq 0$. We can identify
$\F(\al)$ with functions on $\Z$ and want to consider the operator
$L$ with respect to a rescaled basis $\de_l = h_l\, e_l$,
$h_l\in\C$, such that $L$ becomes a symmetric operator when
considered on $\ell^2(\Z)$ with the standard orthonormal basis $\{
e_l\}_{l\in\Z}$. By a straightforward calculation we see that this
requires
\begin{equation}\label{eq:conditiononhls}
\frac{h_{l+1}^2}{h_l^2} = \al^2 q^{2l+1}  \frac{(\be+\al^2
q^{2l})(1+\al^2q^{2l+2})}{(1+\al^2 q^{2l})(1+\be \al^2 q^{2l+2})},
\end{equation}
and then $L=L(\al,\be)$ acts in the basis $\{e_l\}_{l\in\Z}$ by
\begin{equation}\label{eq:defLsymmform}
L\, e_l = a_l e_{l+1} + b_l e_l + a_{l-1}e_{l-1},
\end{equation}
where
\begin{equation}
\label{eq:a_l}
 a_l=\frac{\al q^{l+\frac12}}{1+\al^2 q^{2l+1}}
\sqrt{\frac{(\be + \al^2 q^{2l})(1+\be\al^2 q^{2l+2})}{(1+\al^2
q^{2l})(1+\al^2 q^{2l+2})}}
\end{equation}
and
\begin{equation}
\label{eq:b_l}
b_l=1-\frac{1+\be\al^2 q^{2l}}{(1+\al^2
q^{2l})(1+\al^2 q^{2l-1})} - \frac{1+\be \al^{-2} q^{-2l}}{(1+
\al^{-2}q^{-2l})(1+\al^{-2}q^{-2l-1})}.
\end{equation}
It follows immediately from \eqref{eq:a_l} that
\begin{equation}\label{eq:asymptoticsals}
a_l=\begin{cases} \al \sqrt{\be}q^{l+\frac12} + {\mathcal O}(q^{3l}),
& l\to\infty, \\[3mm]
\al^{-1} \sqrt{\be}q^{-l-\frac12} + {\mathcal O}(q^{-3l}), &l\to-\infty.
      \end{cases}
\end{equation}
Similarly, after simplifying \eqref{eq:b_l} to
\begin{equation}
\label{eq:b_l new} b_l=\frac{\al^2(1+q)(1-\be
q)q^{2l-1}}{(1+\al^2q^{2l+1})(1+\al^2q^{2l-1})},
\end{equation}
it follows that
\begin{equation}\label{eq:asymptoticsbls}
b_l = \begin{cases} \al^2(1+q)(1-\be q) q^{2l-1} + {\mathcal O}(q^{4l}), & l\to\infty,\\[3mm]
\al^{-2}(1+q)(1-\be q) q^{-2l-1} + {\mathcal O}(q^{-4l}), & l\to
-\infty.
       \end{cases}
\end{equation}

Stressing the dependence of the coefficients on $\al$ and $\be$, we
have
\begin{equation}
\label{eq:symmetry1} a_l(\al,1/\be q^{2}) =
\frac{a_l(\al,\be)}{q\be}, \quad b_l(\al,1/\be q^{2}) =
-\frac{b_l(\al,\be)}{q\be},
\end{equation}
and this implies that
\begin{equation}
\label{eq:U}
-q\be \, L(\al,1/\be q^{2}) = U L(\al,\be) U^\ast,
\end{equation}
where $U\in {\mathcal B}\bigl(\ell^2(\Z)\bigl)$ is the unitary
involution given by $Ue_l = (-1)^l e_l$. Moreover, since
$a_l(-\al,\be) = -a_l(\al,\be)$ and $b_l(-\al,\be) = b_l(\al,\be)$,
\begin{equation}
\label{eq:U2}
L(-\al,\be) = U L(\al,\be) U^\ast.
\end{equation}
Furthermore,
$a_{-l-1}({1}/{\al},\be) = a_l(\al,\be)$ and $b_{-l}({1}/{\al},\be) = b_l(\al,\be)$,
which implies
\begin{equation}
\label{eq:V}
L({1}/{\al}, \be)= V L(\al,\be) V^\ast
\end{equation}
for $V\in {\mathcal B}\bigl(\ell^2(\Z)\bigr)$ given by
$Ve_l=e_{-l}$.

As a consequence of \eqref{eq:U}, we can restrict
ourselves to the case $0\leq \be\leq 1/q$ (or $\be\geq 1/q$) without loss
of generality. Note that $b_l(\al,\be)=0$ for all $l\in\Z$ precisely
when $\beta=1/q$.

\begin{prop}\label{prop:Lselfadjointtraceclass}
$L$ defined by \eqref{eq:defLsymmform}--\eqref{eq:b_l} on the dense
subspace ${\mathcal D} \subset \ell^2(\Z)$ of finite linear
combinations of elements from the orthonormal basis $\{
e_l\}_{l\in\Z}$ extends to a bounded self-adjoint operator on
$\ell^2(\Z)$. Moreover, $L$ is a compact operator with singular
values $s_n(L) = {\mathcal O}(q^{n/2})$. In particular, $L$ belongs
to every Schatten $p$-class, $0<p<\infty$.
\end{prop}

For the general theory of operator ideals we refer to Gohberg and
Kre\u{\i}n \cite{GohbK}. However, we find it more convenient to
number the singular values starting from $n=0$.

\begin{proof}

From the explicit forms of $a_l$ and $b_l$ in
\eqref{eq:a_l}--\eqref{eq:b_l} and the estimates in
\eqref{eq:asymptoticsals}, \eqref{eq:asymptoticsbls} it follows that
the sequences $\{ a_l\}_{l\in\Z}$ and $\{b_l\}_{l\in\Z}$ are
bounded, so that $L$ extends to a bounded operator on $\ell^2(\Z)$.
The self-adjointness of $L$ follows directly from the fact that
$\langle L\, e_l,e_k\rangle = \langle  e_l,L\, e_k\rangle$ for all
$l,k\in\Z$, since $a_l, b_l\in\R$.

To prove the second part of the proposition observe that for an
operator $L$ of the form as in \eqref{eq:defLsymmform} the operator
norm satisfies, cf. \cite[Lemma 3.3.3]{Koel} for the case of
Jacobi operators,
\begin{equation}\label{eq:normestimate}
\| L \| \leq 2 \sup_{l\in \Z}|a_l| + 2 \sup_{l\in \Z}|b_l|
\end{equation}
Define $P_{2k+1}$, $k\geq 0$, as the orthogonal projection on
$\text{span}\{ e_{-k}, e_{-k+1},\ldots, e_k\}$ and $P_{2k}$, $k\geq
0$, as the orthogonal projection on $\text{span}\{ e_{-k+1},
e_{-k+2},\ldots, e_{k}\}$ with the convention that $P_0$ is the
zero-projection. Then $\text{dim}\, \text{Ran}(P_n)= n$ for all $n\geq 0$.

First note that, by \eqref{eq:normestimate},
\eqref{eq:asymptoticsals} and \eqref{eq:asymptoticsbls},
\begin{equation*}
\| L\,  -\, P_{2k+1}L\| \leq  2 \sup_{l\in \Z, |l|\geq k}|a_l| +
2 \sup_{l\in \Z, |l|\geq k}|b_l|  = {\mathcal O}(q^k),
\end{equation*}
so that $L$ can be approximated in operator norm by the finite rank
operators $P_{2k+1}L$. This implies that $L$ is a compact operator.
We similarly have
\begin{equation*}
\| L\,  -\, P_{2k}L\| \leq   2 \sup_{l\in \Z, |l|\geq k}|a_l| +
2 \sup_{l\in \Z, |l|\geq k}|b_l| = {\mathcal O}(q^k).
\end{equation*}
By \cite[Thm. 2.1]{GohbK} the $n$-th singular value is
\begin{equation*}
 s_n(L) = \min\{\| L- K\| \mid K\, \text{finite rank operator},
\text{dim}\,\text{Ran}(K)\leq n\}
\end{equation*}
and by taking $K=P_n L$ we get the desired result.
\end{proof}

\subsection{Eigenfunctions for $L$}\label{ssec:eigenfunctionsforL}

In order to derive the spectral decomposition for the self-adjoint
operator $L$ we use the approach of finding the resolvent in terms
of suitably behaved eigenfunctions, see e.g. \cite{Koel}. So we look
for solutions to
\begin{equation}\label{eq:eigenvalueequation}
z\, \psi_l(z)=a_l\psi_{l+1}(z)+b_l\psi_l(z)+a_{l-1}\psi_{l-1}(z),
\quad l\in\Z,
\end{equation}
and, more specifically, for $z\in \C\setminus\R$ we need to
determine the subspaces $S^\pm_z$ consisting of solutions
$\{\psi_l(z)\}_{l\in\Z}$ to \eqref{eq:eigenvalueequation} such that
$\sum_{l=N}^\infty |\psi_{\pm l}(z)|^2<\infty$ for some $N\in \Z$.
Clearly $\dim S^\pm_z \leq 2$, but since $S^+_z\cap S^-_z$ is the
deficiency space of $L$, Proposition
\ref{prop:Lselfadjointtraceclass} implies that $S^+_z\cap S^-_z =
\{0\}$. We will show in Lemma \ref{lem:manysolutions} that both
$S^+_z$ and $S^-_z$ are non-trivial, and it therefore follows that
$\dim S^\pm_z = 1$.

Since the Al-Salam--Chihara polynomials are related to the little
$q$-Jacobi polynomials, cf. the remark in the introduction, one may
expect to be able to find solutions to \eqref{eq:eigenvalueequation}
in terms of associated little $q$-Jacobi polynomials. These
polynomials have been studied by Gupta, Ismail and Masson
\cite{GuptIM} as a limiting case of the associated big $q$-Jacobi
polynomials. Instead of taking the appropriate limit in their
results, we take the limit in the corresponding contiguous
relations, which are attributed to Libis in \cite{GuptIM}. This
leads to contiguous relations for ${}_2\vp_1$-series.

We take $e=cw$ in \cite[(2.5)]{GuptIM}, let $c\to 0$ and relabel to
find the relation
\begin{equation}\label{eq:Libiscontiguousrel1}
\begin{split}
\Psi^- - &\left( \frac{(c-a)(c-b)}{(1-c)(q-c)}z +\,
\frac{(q-a)(q-b)}{(q^2-c)(q-c)}zq + (1-z)\right) \Psi
\\ &\qquad + \frac{(c-a)(c-b)(1-a)(1-b)}{(1-c)^2(1-qc)(q-c)}
z^2 \, \Psi^+ = 0,
\end{split}
\end{equation}
where
\begin{equation*}
\Psi = \tfo{a,b}{c}{z}, \quad \Psi^\pm = \tfo{aq^{\pm 1},bq^{\pm
1}}{cq^{\pm 2}}{z}.
\end{equation*}

Another contiguous relation,
\begin{equation}
\label{eq:GRcontiguousrel1}
\biggl(\frac{q(1+q)(ab-q)}{(aq-b)(a-bq)}+z\biggr)\phi=
\frac{q(1-a)(b+q)}{(a-b)(b-aq)}\,\phi^+ +
\frac{q(1-b)(a+q)}{(b-a)(a-bq)}\,\phi^-,
\end{equation}
where
\begin{equation*}
\phi = \tfo{a,b}{-q}{z}, \quad \phi^+ = \tfo{aq,b/q}{-q}{z}, \quad
\phi^- = \tfo{a/q,bq}{-q}{z},
\end{equation*}
can be obtained by combining identities from \cite[Ex.1.10]{GaspR}
or else be verified directly.

In the following lemma we use the notation
$\theta(z)=(z,q/z;q)_\infty$, $z\neq 0$, for the rescaled Jacobi
$\theta$-function. It satisfies the relation
\begin{equation}\label{eq:thetaproperty}
\theta(zq^l)=(-z)^{-l} q^{-\binom{l}{2}}\theta(z)
\end{equation}
for all $l\in\Z$.

\begin{lemma}\label{lem:manysolutions}
\text{\rm (i)} The functions $\psi_l(z) = \psi_l(z;\al,\be)$ defined
by
\begin{equation}
\label{eq:psi}
\psi_l(z)=C_l(\al,\be)\,\al^l\be^{l/2}q^{l^2/2}\frac{(1/z;q)_\infty}
{z^l} \tfo{i\al\sqrt{\be}q^{l+1},\,
-i\al\sqrt{\be}q^{l+1}}{-\al^2q^{2l+1}}{\frac{1}{z}}
\end{equation}
with
\[
\frac{C_l(\al,\be)}{\sqrt{1+\al^2q^{2l}}}=
\frac{\sqrt{(-\al^2/\be,-\al^2\be q^2;q^2)_l}}{(-\al^2q;q)_{2l}}
\]
satisfy \eqref{eq:eigenvalueequation} for $|z|>1$.
Moreover, each $\psi_l(z)$ has an analytic continuation to
$\C\setminus\{0\}$ such that $\{\psi_l(z)\}_{l\in\Z}\in S^+_z$.

The functions $\Psi_l(z) = \Psi_l(z;\al,\be) :=
\psi_{-l}(z;{1}/{\al},\be)$ also satisfy
\eqref{eq:eigenvalueequation} for $|z|>1$, and each $\Psi_l(z)$ has
an analytic continuation to $\C\setminus\{0\}$ such that
$\{\Psi_l(z)\}_{l\in\Z}\in S^-_z$.

\text{\rm (ii)} The functions $\phi_l(z) =\phi_l(z;\al,\be)$ defined
by
\begin{equation}
\label{eq:phi}
\phi_l(z) = B_l(\al,\be)\, q^{-l/2} \sqrt{1+\al^2q^{2l}}
\, \tfo{i\al q^l/\sqrt{\be}, -iq^{-l}/\al\sqrt{\be}}{-q}{qz},
\end{equation}
where $B_l(\al,\be)$ is a phase factor given by
\[
B_l(\al,\be) = (-i)^l e^{i\arg\left(w_l(\al,\be)\right)}, \quad
w_l(\al,\be)= ({i\al}/{\sqrt{\be}}, -i\al\sqrt{\be}q;q)_l,
\]
satisfy \eqref{eq:eigenvalueequation} for $|z|<1/q$.

\text{\rm (iii)} The functions in \eqref{eq:psi} and \eqref{eq:phi}
are related by
\[
\psi_l(z;\al,\be) = K(z;\al,\be) \phi_l(z;\al,\be) +
K(z;-\al,\be) (-1)^l\phi_l(z;-\al,\be)
\]
with
\[
K(z;\al,\be) =
\frac{(-i\al/\sqrt{\be}, i\al\sqrt{\be}q;q)_\infty}{2(-q,-\al^2 q;q)_\infty}
\frac{\theta(iz/\al\sqrt{\be})}{(-z/\be;q)_\infty}.
\]
When $\be>1$, this gives the analytic continuation of $\psi_l(z)$ to
$\C\setminus\{0\}$. Moreover,
\begin{multline}
\label{eq:all3}
 \qquad \theta(i\al zq/\sqrt{\be})
(-\al^2q,-i/\al\sqrt{\be},i\sqrt{\be}q/\al;q)_\infty \,\psi_l(z) \\+
\theta(i\al\sqrt{\be}/z)
(-1/\al^2,i\al/\sqrt{\be},-i\al\sqrt{\be}q;q)_\infty \,\Psi_l(z)
 = \\
\theta(1/z)\theta(-1/\al^2)(-q,-\be q/z;q)_\infty \,\phi_l(z).
\qquad
\end{multline}
\end{lemma}

We can use the symmetries \eqref{eq:U}, \eqref{eq:U2}, and
\eqref{eq:V} to write down many more solutions, in particular
$(-1)^l\phi_l(z;-\al,\be)$ from Lemma \ref{lem:manysolutions}(iii)
is a solution to \eqref{eq:eigenvalueequation}. However, most of
these solutions can also be obtained from standard transformation
formulae for basic hypergeometric series, such as in \cite[App.
III]{GaspR}.

The function $K(\,\cdot\,;\pm\al,\be)$ has an essential singularity
at $0$ and simple poles at $-\be q^{-k}$, $k\geq 0$. So in case
$\be>1$, Lemma \ref{lem:manysolutions}(iii) gives the analytic
extension of $\psi_l(z)$ to $\C\setminus\{0\}$. But the limit case
$\be\da 0$ is an important special case. For $l$ sufficiently large,
however, the analytic extension of $\psi_l(z)$ to $\C\setminus\{0\}$
follows easily from Heine's transformation \cite[(1.4.1)]{GaspR},
\begin{multline}\label{eq:altanalyticext}
(1/z;q)_\infty \tfo{i\al\sqrt{\be}q^{l+1},\,
-i\al\sqrt{\be}q^{l+1}}{-\al^2q^{2l+1}}{\frac{1}{z}} \\ =
\frac{(-i\al\sqrt{\be}q^{l+1},i\al\sqrt{\be}q^{l+1}/z;q)_\infty}
{(-\al^2q^{2l+1};q)_\infty} \tfo{1/z,\, -i\al q^l/\sqrt{\be}}
{i\al\sqrt{\be}q^{l+1}/z}{-i\al\sqrt{\be}q^{l+1}}
\end{multline}
and the fact that $(c;q)_\infty {}_2\vp_1(a,b;c;q,z)$ is analytic in
its parameters.

Note that the phase factor in Lemma \ref{lem:manysolutions}(ii) also
can be written as
\[
B_l(\al,\be)=(-i)^l \sqrt{\frac{({i\al}/{\sqrt{\be}}, -i\al\sqrt{\be}q;q)_l}
{(-{i\al}/{\sqrt{\be}}, i\al\sqrt{\be}q;q)_l}}.
\]
For future reference we write down an explicit expression for
$\Psi_l(z)=\Psi_l(z;\al,\be)$, viz.
\begin{equation}\label{eq:Psi}
\Psi_l(z)= \frac{(1/z;q)_\infty
z^l}{C_l(\al,\be)\,\al^{l}\be^{l/2}q^{l^2/2}}
\tfo{i\sqrt{\be}q^{1-l}/\al ,\,
-i\sqrt{\be}q^{1-l}/\al}{-q^{1-2l}/\al^2}{\frac{1}{z}}.
\end{equation}

\begin{proof}
  The fact that the functions $\psi_l(z;\al,\be)$ for $|z|>1$ and $\phi_l(z;\al,\be)$
  for $|z|<1/q$ satisfy \eqref{eq:eigenvalueequation} follows directly
  from \eqref{eq:Libiscontiguousrel1} and \eqref{eq:GRcontiguousrel1}.
  Since $\psi_l(z)$ and $\phi_l(z)$ are analytic, their analytic
  extensions satisfy \eqref{eq:eigenvalueequation} as well. Invoking
  \eqref{eq:V} shows that $\Psi_l(z)$ also satisfies
  \eqref{eq:eigenvalueequation}.

  The first relation in (iii) follows from Heine's transformation
  formula \cite[(1.4.3)]{GaspR} and the formula \cite[(4.3.2)]{GaspR}
  for the analytic continuation of a ${}_2\vp_1$-series. The second
  relation follows from \cite[(3.3.5)]{GaspR}, using
  \cite[(1.4.5)]{GaspR} for the second ${}_2\vp_1$-series and
  \eqref{eq:thetaproperty}.

The explicit expression of $\psi_l(z)$ in \eqref{eq:psi} shows
immediately that $\{\psi_l(z)\} \in S^+_z$ for $|z|>1$ due to the
factor $q^{l^2/2}$. For the general case we use
\eqref{eq:altanalyticext} to see that $\psi_l(z)$ behaves like
\[
C(\al,\be)z^{-l}\al^l\be^{l/2}q^{l^2/2} \quad \mbox{as } l\ra\infty,
\]
where $C(\al,\be)=\lim_{l\ra\infty}C_l(\al,\be)$.
Once we have the analytic continuation of $\psi_l(z)$ to
$\C\setminus\{0\}$ for some large value of $l$, we can get the
analytic continuation for other values of $l$ by determining
$\psi_l(z)$ recursively from \eqref{eq:eigenvalueequation}.

The statement for $\Psi_l(z)$ follows similarly.
\end{proof}

\begin{remark}
Note that ${}_2\vp_1$-series similar to the ones in Lemma
\ref{lem:manysolutions} already have occurred in the spectral
decomposition of certain self-adjoint operators. Only the parameter
regimes are somewhat different. This is of great influence on the
particular structure of the spectral decomposition. See
\cite{KoelIM} for the case of the curly exponential
$\mathcal{E}_q(z;t)$ related to a \emph{not} essentially
self-adjoint operator. Another related case concerns little
$q$-Jacobi functions and the Askey--Wilson $q$-Hankel transform as
discussed in \cite{KoelSNATO}.
\end{remark}

\subsection{Spectral decomposition of $L$}\label{ssec:spectraldecomp}

The resolvent of $L$ can be written in terms of solutions in
$S^\pm_z$ and their Wronskian, see \eqref{eq:jk} below. Recall that
the Wronskian of two sequences $\phi = \{ \phi_l\}_{l\in\Z}$ and
$\psi= \{ \psi\}_{l\in\Z}$ is the sequence defined by $[\phi,\psi]_l
= a_l (\phi_{l+1}\psi_l - \phi_l\psi_{l+1})$, $l\in\Z$. It is
well-known that $[\phi,\psi]_l$ is independent of $l$ in case both
$\psi$ and $\phi$ satisfy \eqref{eq:eigenvalueequation}, see e.g.
\cite{Koel} for more information.

\begin{prop}\label{prop:Wronskian}
The Wronskian of $\psi(z)=\{\psi_l(z)\}_{l\in\Z}$ and
$\Psi(z)=\{\Psi_l(z)\}_{l\in\Z}$ is given by
\begin{equation}
\label{eq:Wronskian} [\psi(z),\Psi(z)]=-z(-q\be/z,1/z;q)_\infty,
\quad z\neq 0.
\end{equation}
\end{prop}

The proof of the proposition is based on the following lemma.

\begin{lemma}
\label{lem:convergence} Assume that $\be\leq 1/q$. As $x\da 0$, we
have
\[
\tfo{iq\sqrt{\be}x ,\, -iq\sqrt{\be}x}{-x^2q}{1/z}
\ra\frac{1}{(1/z;q)_\infty}
\]
and
\[
\tfo{iq\sqrt{\be}/x ,\, -iq\sqrt{\be}/x}{-q/x^2}{1/z}
\ra(-\be q/z;q)_\infty.
\]
The convergence is uniform for $z$ in compact subsets of $|z|>1$.
\end{lemma}

\begin{proof}
The ${}_2\vp_1$-series can be written as
\[
\usum\frac{(iq\sqrt{\be}x,-iq\sqrt{\be}x;q)_n}{(-x^2q,q;q)_n}(1/z)^n=
1+\sum_{n=1}^\infty\frac{(1+x^2\be q^2)\cdots(1+x^2\be q^{2n})}
{(1+x^2q)\cdots(1+x^2q^n)}\frac{(1/z)^n}{(q;q)_n},
\]
respectively
\[
\usum\frac{(iq\sqrt{\be}/x,-iq\sqrt{\be}/x;q)_n}{(-q/x^2,q;q)_n}(1/z)^n=
1+\sum_{n=1}^\infty\frac{(x^2+\be q^2)\cdots(x^2+\be q^{2n})}
{(x^2+q)\cdots(x^2+q^n)}\frac{(1/z)^n}{(q;q)_n},
\]
using the fact that $(a,-a;q)_n=(a^2;q^2)_n$. So the termwise
convergence is clear recalling Euler's power series expansions of
the $q$-exponential functions, see \cite[(1.3.15-16)]{GaspR}. When
$\be\leq 1/q$ we have $\be q^j\leq 1$ for all $j\geq 1$, and the
result follows by dominated convergence. The statement on uniform
convergence in $|z|>1$ is straightforward.
\end{proof}

Note that the proof easily can be adapted to the case $\be\leq
q^{-k}$, $k\in\N$ fixed.

\begin{proof}[Proof of Proposition \ref{prop:Wronskian}]
Assuming that $|z|>1$ it follows from Lemma \ref{lem:convergence}
and \eqref{eq:psi}, \eqref{eq:Psi} that, as $l\ra\infty$,
\[
\psi_l(z)\sim C_l(\al,\be){\al^l\be^{l/2}q^{l^2/2}}{z^{-l}}
\quad\mbox{and}\quad
\Psi_l(z)\sim\frac{z^l(1/z,-q\be/z;q)_\infty}{C_l(\al,\be)\al^{l}\be^{l/2}q^{l^2/2}}.
\]
Combining this with the asymptotic behavior of $a_l$ from
\eqref{eq:asymptoticsals} and an easy limit for $C_l(\al,\be)$ as
$l\to\infty$, we get the desired expression for the Wronskian. By
analytic continuation and Lemma \ref{lem:manysolutions},
\eqref{eq:Wronskian} remains valid for $z\in\C\setminus\{0\}$.
\end{proof}

We are now in position to determine the spectrum of the compact
operator $L$ defined in \eqref{eq:defLsymmform}--\eqref{eq:b_l}. Our
considerations will also lead to an explicit orthogonal basis for
$\ell^2(\Z)$ consisting of eigenvectors of $L$.

The spectral measure $E$ for $L$ can be
obtained from the resolvents via the formula
\begin{equation}
\label{eq:resolution} \bigl< E\bigl((a,b)\bigr)v,w\bigr>
=\lim_{\delta\da 0}\lim_{\ep\da 0} \frac{1}{2\pi
i}\int_{a+\delta}^{b-\delta} \bigl< (L-s-i\ep)^{-1}v,w \bigr>-
\bigl< (L-s+i\ep)^{-1}v,w \bigr>\, ds,
\end{equation}
valid for all $v,w\in\ell^2(\Z)$, see e.g. \cite{Koel}. The inner
products appearing in the integral can be written as
\begin{equation}
\label{eq:jk} \bigl< \bigl(L-(s\pm i\ep)\bigr)^{-1}v,w \bigr>=
\sum_{k\leq j}\frac{\psi_j(s\pm i\ep)\Psi_k(s\pm i\ep)} {[\psi(s\pm
i\ep),\Psi(s\pm i\ep)]}
(v_k\overline{w}_j+v_j\overline{w}_k)(1-\tfrac{1}{2}\delta_{j,k})
\end{equation}
and since both $\psi_l(z)$ and $\Psi_l(z)$ are analytic in
$\C\setminus\{0\}$, $E$ is concentrated on the zeros of the
Wronskian $[\psi(z),\Psi(z)]$, excluding $z=0$. This is consistent
with Proposition \ref{prop:Lselfadjointtraceclass}.

\begin{thm}
\label{thm:spectrum} The spectrum of $L$ is given by $\si(L)=-\be
q^\N\cup\{0\}\cup q^{\N_0}$. The accumulation point $0$ is not an
eigenvalue of $L$.
\end{thm}

\begin{proof}
The first statement follows from Proposition \ref{prop:Wronskian}.
To prove that $0$ does not belong to the point spectrum of $L$ we
use Lemma \ref{lem:manysolutions}(ii), which tells us that
\[
f_l=B_l(\al,\be) q^{-l/2} \sqrt{1+\al^2 q^{2l}}, \quad l\in\Z
\]
is a solution to \eqref{eq:eigenvalueequation} for $z=0$. Using the
symmetry \eqref{eq:U2}, we see that
\[
g_l=(-1)^lB_l(-\al,\be) q^{-l/2} \sqrt{1+\al^2 q^{2l}}, \quad l\in\Z
\]
is yet another solution. It is clear that neither
$f=\{f_l\}_{l\in\Z}$ nor $g=\{g_l\}_{l\in\Z}$ is square summable. So
it suffices to show that they are linearly independent. A
straightforward calculation using \eqref{eq:asymptoticsals} shows
that
\[
a_l (f_{l+1} g_l - f_l g_{l+1})
\to -2 i \al\sqrt{\be} \quad\mbox{as } l\to\infty.
\]
Therefore, $[f,g]=-2 i \al \sqrt{\be}$ and the linear independence
follows whenever $\be>0$.
\end{proof}

\begin{remark}\label{rmk:betaiszero}
In the study of $L$ we have sometimes assumed that $\be>0$. Let us
briefly remark on the case $\be=0$, showing that the results in
Proposition \ref{prop:Wronskian} and Theorem \ref{thm:spectrum}
remain valid after taking the limit $\be\da 0$.

The solutions in Lemma \ref{lem:manysolutions}(i) can still be
defined when $\be=0$. Using the fact that
\[
\lim_{\be\da 0} C_l(\al,\be) \be^{l/2} = \frac{\al^l
q^{\binom{l}{2}} \sqrt{1+\al^2 q^{2l}}}{(-\al^2q;q)_{2l}}
=:\frac{C_l(\al)}{\al^lq^{l^2/2}}
\]
and the transformation formula
\[
(z;q)_\infty\, {}_2\vp_1(0,0;c;q,z)= {}_0\vp_1(-;c;q,cz),
\]
which is a limiting case of Heine's transformation formula
\cite[(1.4.3)]{GaspR}, we see that
\begin{equation}
\label{eq:psiforbetaiszero} \psi_l(z;\al,0)=\frac{C_l(\al)}{z^l}
\rphis{0}{1}{-}{-\al^2 q^{2l+1}}{- \frac{\al^2 q^{2l+1}}{z}}
\end{equation}
and
\begin{equation}
\label{eq:Psiforbetaiszero} \Psi_l(z;\al,0)=\frac{z^l}{C_l(\al)}
\rphis{0}{1}{-}{-q^{1-2l}/\al^2}{- \frac{q^{1-2l}}{\al^2z}}.
\end{equation}
These two solutions span $S^\pm_z$, respectively. As regards their
Wronskian we observe that
\[
{}_0\vp_1(-;1/x;q;z/x) \to (z;q)_\infty \quad \mbox{as } x\da 0,
\]
so that the ${}_0\vp_1$-series in \eqref{eq:Psiforbetaiszero} tends
to $(1/z;q)_\infty$ as $l\ra\infty$. Using this and the fact that
\[
a_l(\al,0)\sim\al^2q^{2l+1/2} \quad\mbox{as } l\ra\infty,
\]
we get, in accordance with Proposition \ref{prop:Wronskian}, that
the Wronskian is given by
\[
[\psi(z;\al),\Psi(z;\al)]=-z(1/z;q)_\infty, \quad z\neq 0.
\]
The first part of Theorem \ref{thm:spectrum} thus remains valid, but
the proof that $0$ is not an eigenvalue breaks down since
the Wronskian vanishes as $\be\da 0$. This can be
fixed in the following way, however.
Take $z=0$ and $\be = 0$ in \eqref{eq:eigenvalueequation}
using \eqref{eq:a_l}, \eqref{eq:b_l new} and multiply by
\[
{\sqrt{1+\al^2 q^{2l}} (1+\al^2 q^{2l+1})(1+\al^2 q^{2l-1})}
/{\al^2q^{2l-1/2}}
\]
to get
\begin{equation}
\label{eq:0eigenfunction}
0=\frac{q+\al^2q^{2l}}{\sqrt{1+\al^2q^{2l+2}}}\,\psi_{l+1}+
({q}^{\frac12}+{q}^{-\frac12})\sqrt{1+\al^2q^{2l}}\,\psi_l+
\frac{q^{-1}+\al^2q^{2l}}{\sqrt{1+\al^2q^{2l-2}}}\,\psi_{l-1}, \quad
l\in\Z.
\end{equation}
We are looking for two linearly independent solutions to
\eqref{eq:0eigenfunction} that do not belong to $\ell^2(\Z)$. Setting
$\phi_l={(-1)^l\psi_l}/{\sqrt{1+\al^2q^{2l}}}$,
\eqref{eq:0eigenfunction} reduces to
\[
0=(q+\al^2q^{2l})\phi_{l+1}-
({q}^{\frac12}+{q}^{-\frac12})(1+\al^2q^{2l})\phi_l+
(q^{-1}+\al^2q^{2l})\phi_{l-1}, \quad l\in\Z
\]
and one immediately comes up with the solution $\phi_l=q^{-l/2}$.
Another solution is given by
\[
\tilde{\phi}_l=(1+\al^2q^l)(1-q^l)q^{-3l/2}
\]
and whereas $\phi_l$ never vanishes, we have $\tilde{\phi}_0=0$. The
two solutions, and the corresponding solutions to
\eqref{eq:0eigenfunction}, are therefore linearly independent
and obviously give  non-square summable eigenfunctions.
\end{remark}

\subsection{Orthogonality relations}

Theorem \ref{thm:spectrum} shows that the eigenvalues of $L$ form
two sequences, namely $q^n$, $n\geq 0$, and $-\be q^{n+1}$, $n\geq
0$. In the special case $\be=0$ all the negative eigenvalues
disappear. The corresponding eigenfunctions, $\psi(q^n)$ and
$\psi(-\be q^{n+1})$, form an orthogonal basis for $\ell^2(\Z)$
since $L$ is self-adjoint. Note that we can replace $\psi$ by $\Psi$
as the sequences are proportional. The explicit relations are as
follows.

\begin{lemma}\label{lem:proportional}
For $n\geq 0$, we have
\[
(-1)^n\al^{2n+2}\Psi_l(q^n)=\frac{(-\al^2q;q)_\infty}{(-1/\al^2;q)_\infty}
\frac{(-\be q^2/\al^2;q^2)_\infty}{(-\al^2\be q^2;q^2)_\infty}\,\psi_l(q^n) \\
\]
and
\[
(-1)^n \al^{2n+2}\Psi_l(-\be q^{n+1})=
\frac{(-\al^2q;q)_\infty}{(-1/\al^2;q)_\infty}
\frac{(-1/\al^2\be;q^2)_\infty}{(-\al^2/\be;q^2)_\infty}\,
\psi_l(-\be q^{n+1}).
\]
\end{lemma}

\begin{proof}
If we set $z=q^n$ or $z=-\be q^{n+1}$ in \eqref{eq:all3}, the
right-hand side vanishes. Using \eqref{eq:thetaproperty} and
simplifying gives the result.
\end{proof}

As we will see below the eigenfunctions are closely related to the
symmetric Al-Salam--Chihara polynomials. In particular, the fact
that the eigenfunctions are orthogonal can be translated into
orthogonality relations for the polynomials in
\eqref{eq:defsymmhASCpols}.

Applying \cite[(1.4.1)]{GaspR} to the ${}_2\vp_1$-series in
\eqref{eq:psi} we see that $\psi_l(q^n)$ can be written as, see
\eqref{eq:defsymASCpols} and \eqref{eq:defsymmhASCpols},
\begin{equation*}
\frac{(-\al^2q;q)_\infty}{(-\al^2\be q^2;q^2)_\infty}\psi_l(q^n)=
\sqrt{\frac{(-\al^2/\be;q^2)_l}{(-\al^2\be q^2;q^2)_l}}
\al^l\be^{l/2}q^{l^2/2}\sqrt{1+\al^2 q^{2l}}(-1)^n\al^n
\,h_n^{(\be)} \bigl(x_l(\al)|q\bigr).
\end{equation*}
Similarly, but this time applying \cite[(1.4.5)]{GaspR} to the
${}_2\vp_1$-series in \eqref{eq:psi}, we get
\begin{equation*}
\frac{(-\al^2q;q)_\infty}{(-\al^2/\be;q^2)_\infty} \psi_l(-\be
q^{n+1})=\sqrt{\frac{(-\al^2\be q^2;q^2)_l}{(-\al^2/\be;q^2)_l}}
\frac{(-1)^l\al^lq^{l^2/2}}{\be^{l/2}q^l}\sqrt{1+\al^2q^{2l}}
(-1)^n\al^n\,h_n^{(1/\be q^2)}\bigl(x_l(\al)|q\bigr).
\end{equation*}
The fact that $\left<\psi(q^n),\psi(-\be q^{m+1})\right>=0$ for
$n,m\geq 0$ therefore reads
\begin{equation}\label{eq:mixedOR}
\sum_{l=-\infty}^\infty
(-1)^l\alpha^{2l}(1+\alpha^2q^{2l})q^{l(l-1)}
h_n^{(\beta)}(x_l(\alpha)|q)h_m^{(1/\beta q^2)}(x_l(\alpha)|q)=0.
\end{equation}
A direct proof of \eqref{eq:mixedOR} can also be given since
\begin{equation*}
\sum_{l=-\infty}^\infty
(-1)^l\alpha^{2l}(1+\alpha^2q^{2l})q^{l(l-1)} x_l^k(\alpha)=0,
\qquad k\geq 0,
\end{equation*}
using Jacobi's triple product \cite[(1.6.1)]{GaspR} and the binomial
expansion.


In order to write down the relation
$\left<\psi(q^n),\psi(q^m)\right>=\de_{n,m}\Vert\psi(q^n)\Vert^2$
explicitly, we need to determine the $\ell^2(\Z)$-norm of the
sequence $\psi(q^n)=\{\psi_l(q^n)\}_{l\in\Z}$. This will be an easy
task once we know the spectral projection $E\bigl(\{q^n\}\bigr)$
given by, cf. \eqref{eq:resolution},
\[
\left<E\bigl(\{q^n\}\bigr)v,w\right>=-\frac{1}{2\pi i}\oint_{(q^n)}
\bigl<(L-s)^{-1}v,w \bigr>\,ds, \quad v,w\in\ell^2(\Z).
\]
Alluding to \eqref{eq:jk} and the fact that
\[
\frac{1}{2\pi i}\oint_{(q^n)}
\frac{\psi_j(s)\Psi_k(s)}{[\psi(s),\Psi(s)]}\,ds=
\psi_j(q^n)\Psi_k(q^n)\,\,\underset{z=q^n}{\text{Res}}\,\frac{1}{[\psi(z),\Psi(z)]},
\]
it follows from Lemma \ref{lem:proportional} that
\begin{equation}
\label{eq:spectralproj}
\left< E\bigl(\{q^n\}\bigr)v,w \right>=
\frac{q^{n^2} \bigl<v,\psi(q^n)\bigr>\,\bigl<\psi(q^n),w\bigr>}
{\al^{2n+2}\be^n(-1/\be,q;q)_n} \frac{(-\al^2q;q)_\infty(-\be
q^2/\al^2;q^2)_\infty} {(-1/\al^2,-\be q,q;q)_\infty(-\al^2\be
q^2;q^2)_\infty}
\end{equation}
since
\[
\underset{z=q^n}{\text{Res}}\frac{1}{[\psi(z),\Psi(z)]}
=\frac{(-1)^{n+1}q^{n^2}}{\be^n(-1/\be,q;q)_n}\frac{1}{(-\be
q,q;q)_\infty}.
\]
In particular, setting $v=w=\psi(q^n)$ we see that
\[
\Vert\psi(q^n)\Vert^2= \frac{\al^{2n+2}\be^n(-1/\be,q;q)_n}{q^{n^2}}
\frac{(-1/\al^2,-\be q,q;q)_\infty(-\al^2\be
q^2;q^2)_\infty}{(-\al^2q;q)_\infty(-\be q^2/\al^2;q^2)_\infty}
\]
and the orthogonality relation $\langle \psi(q^n),\psi(q^m)\rangle =
\de_{n,m} \| \psi(q^n)\|^2$ reads
\begin{multline}
\label{eq:OR} \qquad
\sum_{l=-\infty}^\infty\frac{(-\al^2/\be;q^2)_l}{(-\al^2\be
q^2;q^2)_l} \al^{2l}\be^l(1+\al^2q^{2l})q^{l^2}
h_n^{(\be)}(x_l(\al)|q)h_m^{(\be)}(x_l(\al)|q)\\
=\de_{n,m}\frac{\be^n (-1/\be,q;q)_n}{q^{n^2}}
\frac{(-\al^2,-q/\al^2,-\be q,q;q)_\infty}{(-\al^2\be q^2,-\be
q^2/\al^2;q^2)_\infty}. \qquad
\end{multline}
It follows from \eqref{eq:spectralproj} that the spectral projection
$E\bigl(\{q^n\}\bigr)$ has rank $1$. In other words, the eigenvalues
$q^n$, $n\geq 0$, are simple.

A similar computation can be carried out to determine the
$\ell^2(\Z)$-norm of $\psi(-\be q^{n+1})$, but the relation
$\left<\psi(-\be q^{n+1}),\psi(-\be
q^{m+1})\right>=\de_{n,m}\Vert\psi(-\be q^{n+1})\Vert^2$ corresponds
to \eqref{eq:OR} with $\be$ replaced by $1/\be q^2$.

The orthogonality relation \eqref{eq:OR} was obtained by
Christiansen and Ismail in \cite[Thm. 6.1]{ChriI} using the
attachment procedure combining moments with generating functions and
mimicking the construction of the $N$-extremal solutions to the
$q^{-1}$-Hermite moment problem. The present analysis, however,
gives a much more complete picture of the situation. In particular,
we obtain an explicit orthogonal basis for $L^2(\la_\al^{(\be)})$,
where $\la_\al^{(\be)}$ denotes the discrete probability measure in
\eqref{eq:OR}.

\begin{thm}\label{thm:basisforweightedL2}
The polynomials are dense in $L^2(\la_\al^{(\be)})$ if and only if
$\beta=0$. When $\beta=0$, the polynomials $\{h_n(x|q)\}_{n\geq 0}$
form an orthogonal basis. For $\beta>0$, let $\Phi^{(\beta)}(x|q)$
denote the function given by
\[
\Phi^{(\beta)}(\sinh y|q)=
\frac{(-e^{-2y}/\beta;q^2)_\infty}{(-e^{-2y}\beta q^2;q^2)_\infty}
e^{-y\bigl(\frac{i\pi-\log\beta}{\log q}-1\bigr)}, \quad y\in\R.
\]
Then $\{h_n^{(\beta)}(x|q)\}_{n\geq
0}\cup\{\Phi^{(\beta)}(x|q)h_n^{(1/\beta q^2)}(x|q)\}_{n\geq 0}$
form an orthogonal basis for $L^2(\la_\al^{(\be)})$.
\end{thm}

Note that $|\Phi^{(1/\be q^2)}|= |1/\Phi^{(\be)}|$ for $\be>0$.

\begin{proof}
Rewrite \eqref{eq:mixedOR} in terms of
the orthogonality measure $\la_\al^{(\be)}$ to see
that
$$
{(-1)^l}{\be^{-l}q^{-l}} \frac{(-\al^2\be
q^2;q^2)_l}{(-\al^2/\be;q^2)_l} h^{(1/\be q^2)}_m(x_l(\al)|q)
$$
are orthogonal to the symmetric Al-Salam--Chihara polynomials.
Moreover, by \eqref{eq:OR} with $\be$ replaced by $1/q^2\be$, we see
that these functions are orthogonal with respect to the measure
$\la^{(\be)}_\al$. Then combine the fact that
\[
\frac{(-\alpha^2\beta q^2;q^2)_l}{(-\alpha^2/\beta;q^2)_l}
=\frac{(-\alpha^2\beta
q^2;q^2)_\infty}{(-\alpha^2/\beta;q^2)_\infty}
\frac{(-\alpha^2q^{2l}/\beta;q^2)_\infty} {(-\alpha^2\beta
q^{2l+2};q^2)_\infty}
\]
with the identity
\[
q^{l\bigl(\frac{i\pi-\log\beta}{\log q}-1\bigr)}={(-1)^l}{\be^{-l}
q^{-l}}
\]
and recall that $x_l(\al)$ can be written as $\sinh y$ with $e^y=1/\al q^l$.
\end{proof}

Recalling that a solution $\mu$ to an indeterminate moment problem
is $N$-extremal if and only if the polynomials $\C[x]$ are dense in
$L^2(\R,\mu)$, see e.g. \cite[Chapter 2]{Akhi}, we have the
following consequence of the above theorem.

\begin{cor}\label{cor:nonNextremalASC}
For $\be>0$ the measures $\la_\al^{(\be)}$, $\al\in(q,1]$, are
mutually different, non $N$-extremal solutions to the moment problem
associated with the symmetric Al-Salam--Chihara polynomials.
\end{cor}

Christiansen and Ismail \cite[\S 7]{ChriI} proved this result for
$\al=\be=1$ in a different way. Note that our approach even tells that
the measure $\la_\al^{(\be)}$ cannot be canonical of any order, see
e.g. \cite{Akhi} for more information on canonical measures. This
follows from the fact that the closure of the polynomials has
codimension $+\infty$ in $L^2(\la_\al^{(\be)})$.

\begin{remark}\label{rmk:othercases}
  In \cite{ChriPhD} a classification of the indeterminate moment
  problems within the $q$-analogue of the Askey scheme is given. Since
  each of the corresponding orthogonal polynomials also satisfies a
  second order difference equation, it is natural to look for a
  corresponding self-adjoint operator. This is done for the
  Stieltjes--Wigert polynomials in \cite{ChriK}, and indirectly also
  in \cite{CiccKK} for the $q$-Laguerre polynomials.  It seems that
  for the Al-Salam--Chihara polynomials as in \eqref{eq:defASCpols}
  we need at least $a+b\in i\R$ and $ab\in\R$, which are fulfilled for
  $a=-b\in\R$.  However, in general it seems difficult to find a
  suitable Hilbert space of functions on which the second order
  difference operator can be realized as a self-adjoint operator, and
  some external guidance for such a choice seems necessary. See also
  \cite{KoelS} for an example involving a different type of Hilbert
  space related to the continuous dual $q^{-1}$-Hahn polynomials, and
  more generally \cite{KoelSNATO} for an overview.
\end{remark}

\begin{remark}\label{rmk:dualorthogonality}
Since Theorem \ref{thm:basisforweightedL2} gives an explicit
orthogonal basis for the space $L^2(\la^{(\be)}_\al)$, we also
have the dual orthogonality relations
\begin{equation}
\sum_{n=0}^\infty \frac{\psi_l(q^n) \psi_k(q^n)}{\|\psi(q^n)\|^2} +
\sum_{n=0}^\infty \frac{\psi_l(-\be q^{n+1}) \psi_k(-\be
q^{n+1})}{\|\psi(-\be q^{n+1})\|^2} = \de_{k,l} \quad\text{for }
k,l\in\Z.
\end{equation}
Each of the series can be summed as a special case of the bilinear
generating function \cite[Thm.~7.2]{ChriI}, and we find an identity
involving four balanced ${}_4\vp_3$-series. We do not work out the
details.
\end{remark}

\begin{remark}\label{rmk:directintegral}
The operator $L$ can also be considered as an operator acting on
$L^2(\R)$, and this case can be reduced to a direct integral of the
cases studied in Theorems \ref{thm:spectrum} and
\ref{thm:basisforweightedL2} since $\cup_{\al \in(q,1]}Z(\al)=\R$.
Because the spectrum of $L(\al,\be)$ is independent of $\al$, this
gives no new results. We merely see that the polynomials
$h_n^{(\be)}(x|q)$ are orthogonal with respect to the weight
function
\[
w(x)=\frac{1}{(-e^{2y}/\be,-e^{-2y}/\be;q^2)_\infty}, \quad x=\sinh
y\in\R,
\]
a result contained in \cite[Thm. 5.1]{ChriI}. Observe namely that
\[
w\bigl(x_l(\al)\bigr)= \frac{(-\al^2/\be;q^2)_l}{(-\al^2\be
q^2;q^2)_l}
\frac{\al^{2l}\be^lq^{l(l+1)}}{(-1/\al^2\be,-\al^2/\be;q^2)_\infty}
\]
so that
\[
\lambda_\al^{(\be)}\bigl(\{x_l(\al)\}\bigr)=
M(\al,\be)\,w\bigl(x_l(\al)\bigr)\sqrt{x_l^2(\al)+1}
\]
for some constant $M(\al,\be)$ not depending on $l$. The factor
$\sqrt{x_l^2(\al)+1}$ comes from the change of variables $x=\sinh
y$.

See \cite{ChriK} for direct integral techniques applied to the case
of the Stieltjes--Wigert polynomials.
\end{remark}


\section{Special cases and additional results}

\subsection*{Continuous $q^{-1}$-Hermite polynomials}

Theorem \ref{thm:basisforweightedL2} in case $\be=0$ is dealing with
the continuous $q^{-1}$-Hermite polynomials,
$h_n^{(0)}(x|q)=h_n(x|q)$. These polynomials were introduced by Askey
\cite{Aske} and the associated indeterminate moment problem was
studied by Ismail and Masson \cite{IsmaM} in much detail. In
particular, Ismail and Masson were able to find all of the
$N$-extremal solutions explicitly using the Nevanlinna
parametrization. Corollary \ref{cor:NextremalcontqHermite} below is
thus a result due to Ismail and Masson \cite[\S 6]{IsmaM}, but we give
a different proof based on Theorem \ref{thm:basisforweightedL2}.

\begin{cor}\label{cor:NextremalcontqHermite}
The measures $\la_\al^{(0)}$, $\al\in(q,1]$, constitute all of the
$N$-extremal solutions to the $q^{-1}$-Hermite moment problem.
\end{cor}

\begin{proof}
Theorem \ref{thm:basisforweightedL2} shows that the polynomials are
dense in $L^2(\la_\al^{(0)})$ and since this characterizes the
$N$-extremal measures, it follows that each $\la_\al^{(0)}$ is
$N$-extremal. Now $N$-extremal measures are discrete and for each
$x\in\R$ there is a unique $N$-extremal measure, say $\rho$, with
$\rho(\{x\})>0$, see e.g. \cite{Akhi}. Given $x\in\R$, pick
$\al\in(q,1]$ such that $x=x_l(\al)$ for some $l\in\Z$. Then
$\la_\al^{(0)}(\{x\})>0$ and we have obtained all of the
$N$-extremal measures.
\end{proof}

Note that the explicit orthogonality relation for the
continuous $q^{-1}$-Hermite polynomials with respect to
the measure $\la^{(0)}_\al$ is
\begin{equation}\label{eq:orthocontqinvHermite}
\sum_{l=-\infty}^\infty \al^{4l}(1+\al^2 q^{2l})q^{l(2l-1)}
\, h_n(x_l(\al)|q) h_m(x_l(\al)|q) = \de_{n,m}
\frac{(q;q)_n}{q^{\binom{n+1}{2}}} (-\al^2,-q/\al^2,q;q)_\infty.
\end{equation}

In the special case $\be=1/q$, the fixed point under the involution
$\be\mapsto 1/q^2\be$, we are also dealing with the continuous
$q^{-1}$-Hermite polynomials since $h^{(1/q)}_n(x|q)=h_n(x|q^2)$ as
observed in Section \ref{sec:ASCpolynomials}. In this case the
eigenvectors corresponding to the eigenvalues $q^n$ and $-q^n$ differ
only by a sign,
\begin{equation}
\label{eq:casebeis1/q}
\psi_l(q^n;\al,1/q)=
\al^lq^{\binom{l}{2}}\sqrt{1+\al^2 q^{2l}}(-1)^n\al^n h_n^{(1/q)}(x_l(\al)|q)=
(-1)^l \psi_l(-q^n;\al,{1}/{q}).
\end{equation}
The orthogonality relations \eqref{eq:OR} and \eqref{eq:mixedOR} then take
the simpler forms
\begin{multline}
\label{eq:ORcasebeis1/q}
\qquad\qquad
\sum_{l=-\infty}^\infty \al^{2l}(1+\al^2q^{2l})q^{l(l-1)}
h_n^{(1/q)}(x_l(\al)|q)h_m^{(1/q)}(x_l(\al)|q) \\
=2\de_{n,m}\frac{(q^2;q^2)_n}{q^{n(n+1)}}(-\al^2,-q^2/\al^2,q^2;q^2)_\infty
\qquad\qquad
\end{multline}
and
\begin{equation}
\label{eq:mixedORcasebeis1/q}
\sum_{l=-\infty}^\infty (-1)^l\al^{2l}(1+\al^2q^{2l})q^{l(l-1)}
h_n^{(1/q)}(x_l(\al)|q)h_m^{(1/q)}(x_l(\al)|q) = 0.
\end{equation}
Adding \eqref{eq:ORcasebeis1/q} and \eqref{eq:mixedORcasebeis1/q}
leaves only the even terms in the sum, and since $x_{2l}(\al)$ in base
$q$ equals $x_l(\al)$ in base $q^2$, we obtain
\eqref{eq:orthocontqinvHermite} in base $q^2$. Similarly, subtracting
\eqref{eq:mixedORcasebeis1/q} from \eqref{eq:ORcasebeis1/q} leaves
only the odd terms and we find \eqref{eq:orthocontqinvHermite} in base
$q^2$ with $\al$ replaced by $\al q$.

This can be interpreted in the operator theoretic approach as follows.
We have $b_l(\al,1/q)=0$ and it thus follows that $L^2$ leaves the
subspaces $M^e$ spanned by $\{e_{2l}\}_{l\in\Z}$ and $M^o$ spanned by
$\{e_{2l+1}\}_{l\in\Z}$ invariant. Stressing the dependence on $\al$,
$\be$, and $q$, it is an easy verification that
\[
L(\al,{1}/{q}|q)^2\big\vert_{M^e} \simeq L(\al,0|q^2)
\quad\mbox{and}\quad
L(\al,{1}/{q}|q)^2\big\vert_{M^o} \simeq
L(\al q,0|q^2),
\]
using the fact that $M^e\simeq \ell^2(\Z) \simeq M^o$ as subspaces of
$\ell^2(\Z)$. Therefore,
\[
L(\al,1/q|q)^2 \simeq L(\al,0|q^2) \oplus
L(\al q,0|q^2).
\]
By Theorem \ref{thm:spectrum} the point spectrum of $L(\al,1/q|q)^2$
equals $q^{2\N_0}$, each point having multiplicity $2$. Also, by
Theorem \ref{thm:spectrum}, $L(\al,0|q^2)$ and $L(\al q,0|q^2)$ have
spectrum $q^{2\N_0}$, with each point of multiplicity $1$. Since we
have the eigenvectors explicitly, and because of
\eqref{eq:casebeis1/q}, we obtain
\[
\psi(q^n;\al,1/q|q) + \psi(-q^n;\al,1/q|q) =
2(-\al^2 q^2;q^2)_\infty\, \psi(q^{2n};\al,0|q^2)
\]
and
\[
\psi(q^n;\al,1/q|q) - \psi(-q^n;\al,1/q|q)=
2\al q^{-n} (-\al^2q^4;q^2)_\infty\, \psi(q^{2n};\al q,0|q^2),
\]
stressing the dependence on the base as well. So the measure
$\la^{(1/q)}_\al$ splits as a convex linear combination of two
$N$-extremal measures, as observed by Christiansen and Ismail \cite[\S
6]{ChriI}. In particular, the measure in question is \emph{not}
$N$-extremal, and Theorem \ref{thm:basisforweightedL2} and its
preceding proof shows that the orthocomplement of the polynomials is
spanned by the orthogonal functions
\[
h_n(x|q^2)e^{-i\pi y/\log q}, \quad n\geq 0,
\]
where $x=\sinh y$.

\subsection*{Summation formulae}
Recall the generating function, see \cite{ChihI}, \cite{ChriI},
\begin{equation}\label{eq:genfun}
\sum_{n=0}^\infty h^{(\be)}_n(x|q) \frac{q^{\binom{n}{2}}t^n}{(q;q)_n}
= \frac{(te^{-y},-te^{-y};q^2)_\infty}{(-t^2\be;q)_\infty},
\quad |t|<\frac{1}{\sqrt{\be}}.
\end{equation}
This and Bailey's ${}_6\psi_6$-summation formula, see
\cite[(5.3.1)]{GaspR}, was used by Christiansen and Ismail \cite[\S
6]{ChriI} to obtain the measures $\la_\al^{(\be)}$ as solutions to
the moment problem associated with the symmetric Al-Salam--Chihara
polynomials. Now that we have obtained the measures
$\la_\al^{(\be)}$ as such solutions in a different way, we can
reverse the line of reasoning and obtain a $4$-parameter subfamily
of Bailey's ${}_6\psi_6$-summation formula from \eqref{eq:genfun}.
Multiply \eqref{eq:OR} by $q^{{n(n-1)/2}}t_1^n/(q;q)_n$ and
$q^{{m(m-1)/2}}t_2^m/(q;q)_m$ and sum over $n$ and $m$ using
\eqref{eq:genfun} twice. Interchange the order of summation and
simplify to get
\begin{multline}\label{eq:subBailey}
\,_6\psi_6 \left( \genfrac{.}{.}{0pt}{}
{i\al q,-i\al q,i\al/\sqrt{\be},-i\al/\sqrt{\be},-\al q/t_1,-\al q/t_2}
{i\al,-i\al,i\al\sqrt{\be}q,-i\al\sqrt{\be}q,t_1\al,t_2\al}
\ ;q, \frac{t_1t_2\be}{q} \right) = \\
\frac
{(it_1\sqrt{\be},-it_1\sqrt{\be},it_2\sqrt{\be},-it_2\sqrt{\be},
-t_1t_2/q,-\al^2q,-q/\al^2,-\be q,q;q)_\infty}
{(t_1\al,t_2\al,-t_1/\al,-t_2/\al,t_1t_2\be/q,
i\al\sqrt{\be}q,-i\al\sqrt{\be}q,i\sqrt{\be}q/\al,-i\sqrt{\be}q/\al;q)_\infty}\,,
\qquad
\end{multline}
which is a special case of Bailey's ${}_6\psi_6$-summation
formula that involves $5$ degrees of freedom.

In case we use \eqref{eq:mixedOR} instead of \eqref{eq:OR}, we get the
summation formula
\begin{equation}\label{eq:mixedsubBailey}
\,_4\psi_4 \left( \genfrac{.}{.}{0pt}{}
{i\al q,-i\al q, -\al q/t_1, -\al q/t_2}
{i\al, -i\al, \al t_1, \al t_2}
\ ;q, -\frac{t_1t_2}{q^2} \right) = 0.
\end{equation}
This is the special case $b=qa^{1/2}$ of the ${}_4\psi_4$-summation
formula \cite[(5.3.3)]{GaspR}, which is a special case of Bailey's
${}_6\psi_6$-sum.

We note that one can extend this method by considering suitable
bilinear generating functions as \cite[Thm. 7.2]{ChriI}.

\end{document}